\documentclass[a4paper]{article}

\usepackage{amsmath, amssymb, amsfonts, amscd, amsthm}

\newcommand{\<}{\langle}
\renewcommand{\>}{\rangle}

\renewcommand{\l}{\lambda}
\renewcommand{\L}{\mathbf{\Lambda}}
\newcommand{\B}{\mathbf{B}}
\renewcommand{\b}{\beta}

\renewcommand{\P}{\Psi}

\newcommand{\Z}{\mathbb{Z}}
\newcommand{\C}{\mathbb{C}}

\newcommand{\Q}{\mathbb{Q}}

\newcommand{\F}{\mathbb{F}}

\newcommand{\zbs}{\Z\times BS_\infty^+}
\newcommand{\bs}{BS_\infty^+}
\newcommand{\kba}{\Z\times BGL_\infty(A)^+}
\newcommand{\ba}{BGL_\infty(A)^+}

\newcommand{\spec}{\textrm{Spec}}

\newtheoremstyle{pedro}{}{}{\itshape}{}{\sc}{~--}{ }{\thmname{#1}\thmnumber{ #2}\thmnote{ (#3)}}

\newtheoremstyle{pedrodef}{}{}{\rm}{}{\sc}{~--}{ }{\thmname{#1}\thmnumber{ #2}\thmnote{ (#3)}}

\theoremstyle{pedro}

\newtheorem{lem}{Lemma}[section]

\newtheorem{thm}[lem]{Theorem}

\newtheorem{prop}[lem]{Proposition}

\newtheorem{coro}[lem]{Corollary}

\theoremstyle{remark}

\newtheorem{rmk}[lem]{Remark}

\theoremstyle{pedrodef}

\newtheorem{ex}[lem]{Example}

\newtheorem*{notations2}{Notations}
\newtheorem{defn}[lem]{Definition}

\usepackage{titlesec}

\titleformat{\section}{\bf\center}{\S \arabic{section}.}{1em}{}
\titleformat{\subsection}{\bf}{}{0pt}{}   

\title{Adams operations in cohomotopy}

\author{Pierre Guillot}

\begin{document}


\maketitle

\begin{abstract} We study a collection of operations on the cohomotopy $\pi^0(X)$ of the space $X$, with which it becomes a "$\b$-ring", an algebraic structure analogous to a $\l$-ring. In particular, this ring possesses Adams operations, represented by maps $\bs \to \bs$. We compute their effect in homotopy on the image of $J$, and in mod $2$ cohomology.

The motivation comes from the interpretation of the symmetric group as the general linear group of the "field with one element", which leads to an analogy between cohomotopy and algebraic $K$-theory.

A good deal of this article may be considered as a survey of the theory of $\b$-rings.

\end{abstract}

\section{Introduction}

The symmetric groups $S_n$ and the general linear groups $GL_n(k)$, where $k$ is a field, have a number of features in common. This is particularly clear in the construction of their representations, for example.

A heuristic explanation of this phenomenon involves the "field with one element" $\F_1$ -- see for example the work of Soul\'e \cite{soule}. Let us only say here that one can define a geometry over $\F_1$, which yields combinatorial formulae appearing to be the limit, as $q$ goes to $1$, of classical expressions in the usual geometry over $\F_q$ (eg formulae for the number of rational points). There are several models for this geometry (an alternative to {\em loc cit} being the work of Toen and Vaqui\'e \cite{toen}), though they all agree on one thing: a finite-dimensional vector space over $\F_1$ is simply a finite set (possibly pointed). Thus the symmetric group $S_n$ appears as $GL_n(\F_1)$.

Following this metaphor, we are led to define the algebraic $K$-theory of the field with one element to be the collection of groups $\pi_i(\zbs)$, for $i\ge 0$. By a classical theorem due to Priddy, Quillen, Segal and others (\cite{priddy}), these are the stable homotopy groups of spheres. Thus one is tempted to look for properties in cohomotopy which might resemble those of algebraic $K$-theory.

To illustrate the consistency of all this, let us recall a few classical results. The complex $K$-theory $K(BG)$ of the classifying space of a finite group is the completion of the representation ring $R(G)$ at the augmentation ideal, as was proved by Atiyah \cite{atiyah}. Rector \cite{rector} has proved a similar result involving $K_{\F_q}(BG)$, the $K$-theory of $BG$ over $\F_q$, and the completed ring of modular representations. Finally, "letting $q$ go to $1$", we reach the Segal conjecture, now a theorem by Carlson \cite{carl}, which states that $\pi^0(BG)=[BG,\zbs]$ is the completion of the Burnside ring $A(G)$ at the augmentation ideal.

This has led us to attempt to tackle another family of results. Namely, it is classical that $K(X)$, the complex $K$-theory of the space $X$, is a $\l$-ring (we shall recall the definition below). There is a similar result for $K_A(X)$ ($K$-theory over the ring $A$), see \cite{kra}. We shall prove that $\pi^0(X)$ is a $\b$-ring, which is something similar to a $\l$-ring (but definitely different), following a suggestion of Clemens Berger. For the purposes of this introduction, a $\b$-ring will be a ring $R$ together with operatons $\b_H: R\to R$ indexed by the (conjugacy classes of) subgroups of $S_n$, for various $n$, and satisfying a number of conditions such as $\b_H(\b_K(x))=\b_{H\wr K}(x)$. 

Spending some time investigating the literature on $\b$-rings will lead you to two main observations. First of all, it is surprisingly difficult to unearth the foundations of the theory, with competing definitions appearing in different papers, unfortunate mistakes here which are (sometimes but not always) corrected there, etc. And second, the result that $\pi^0(X)$ is a $\b$-ring is almost already there, though it has never been stated explicitly: all the technical tools needed are available, and strong connections between $\b$-rings and cohomotopy have been established (see \cite{vallejopizero}),

This is why this article may be considered partially as a survey of the theory of $\b$-rings, or perhaps as a clarified, self-contained exposition. We shall strive to make the parallel between $\l$-rings and $\b$-rings as transparent as possible, and this will lead us to a definition which is stronger than that in the literature; however all known examples of $\b$-rings satisfy our extra conditions. 

Having reached this definition, our results are as follows.
\begin{enumerate}
\item $\pi^0(X)$ is a $\b$-ring (theorem \ref{thm:pizero}). Thus there are maps $\b_H:\zbs \to \zbs$ representing the operations.
\item Partial information on the effect of some of these maps on homotopy groups (theorem \ref{thm:imj}). More precisely, we shall consider the so-called Adams operations and calculate their effect on the image of $J$.
\item Partial information on the effect of the Adams operations in cohomology (theorem \ref{thm:psicoho}). It will suffice to say here that, when interpreted correctly, this theorem states that the situation is precisely analogous to that of the Adams operations in complex $K$-theory.
\end{enumerate}

The organization of the paper is rather straightforward, introducing $\l$-rings first, then $\b$-rings, and then proving the main theorem on cohomotopy; after this we explain what Adams operations are before studying their behaviour in cohomology.

\smallskip
\noindent{\em Acknowledgements.} This work was prompted by suggestions from Clemens Berger during my stay at the University of Nice in 2005-2006. I would like to thank him warmly for his enthusiasm and support.

\section{$\l$-rings}

In this section we present the theory of $\l$-rings in a fashion that will make the analogy with $\b$-rings, in the sequel, more transparent. There are plenty of classical expositions in literature: we recomment \cite{atiyahtall} for its concision, \cite{knut} for its clarity, and \cite{patras} for the numerous references it provides.

\subsection{The ring $\L$.}

We define $$\L=\bigoplus_{n\ge 1} R(S_n)$$ where $S_n$ denotes the symmetric group on $n$ letters, and $R(S_n)$ is its representation ring. We ignore the existing multiplication on each $R(S_n)$, and put a ring structure on $\L$ as follows. For $\rho\in R(S_p)$ and $\sigma\in R(S_q)$, we put $$\rho\sigma=Ind_{S_p\times S_q}^{S_{p+q}} (\rho\otimes\sigma).$$ Extending by bilinearity, this endows $\L$ with the structure of a commutative ring with unit, as ones checks easily.

We also have a diagonal map $\Delta:\L\to \L\otimes \L$, defined using the restriction maps $R(S_n) \to R(S_p\times S_q)=R(S_p)\otimes R(S_q)$ (for $p+q=n$). 

The following is well-known (\cite{knut}).

\begin{thm}\label{thm:lambdapoly} Let $\l^n$ denote the signature representation of $S_n$. Then $\L$ is a polynomial ring in the $\l^n$'s: $$\L=\Z[\l^1,\l^2,\ldots,\l^n,\ldots].$$ Similarly, if $\b^n$ is the trivial representation of $S_n$, then $$\L=\Z[\b^1,\b^2,\ldots,\b^n,\ldots].$$
\end{thm}

A direct computation then shows:

\begin{lem}\label{lem:diago} The diagonal map satisfies $$\Delta(\l^n)=\sum_{p+q=n} \l^p\otimes \l^q$$ and
$$\Delta(\b^n)=\sum_{p+q=n} \b^p\otimes \b^q.$$
\end{lem}

In fact, since we will show below that the diagonal map is a ring homomorphism, this lemma determines $\Delta$ completely.

\subsection{The definition of a pre-$\l$-ring.}

If $R$ is any ring, we write $R^R$ for the ring of all maps {\em of sets} $R\to R$.

\begin{defn} A {\em pre-$\l$-ring} is a ring $R$ with a given ring homomorphism $\theta:\L\to R^R$ making the following diagram commute:
$$\begin{CD}
    \L @>{\Delta}>> \L\otimes\L
\\ @V{\theta}VV        @VV{\theta^2}V
\\  R^R @>a>> R^{R\times R} 
\end{CD}$$
\end{defn}
A word of explanation is in order. The map named $\theta^2$ is the only map satisfying, for $a, b\in \L$ and $x, y \in R$, the relation $\theta^2(a\otimes b)(x,y)=\theta(a)(x)\cdot\theta(b)(y)$ (multiplication in $R$). The map $a$ sends $f: R\to R$ to $a(f):R\times R \to R$ with $a(f)(x,y)=f(x+y)$ (the letter $a$ is for addition).

Given theorem \ref{thm:lambdapoly}, we see that in order to define the map $\theta$, we need only specify the operations $\theta(\l^n)$, which are usually written simply $\l^n: R\to R$. Moreover, lemma \ref{lem:diago} asserts that, in a pre-$\l$-ring, one has $$\l^n(x+y)=\sum_{p+q=n}\l^p(x)\cdot \l^q(y)$$ and $$\b^n(x+y)=\sum_{p+q=n}\b^p(x)\cdot \b^q(y).$$ Conversely, when we know that $\Delta$ is a homomorphism of rings, we will deduce that any of these two formulae is sufficient to imply the commutativity of the diagram above.

Thus we reach a very simple definition of a pre-$\l$-ring as a ring with operations $\l^n$ defined for all integers $n$, satisfying the formula above for $\l^n(x+y)$; or equivalently with the operations $\b^n$. However, when we try to parallel the construction for $\b$-rings we shall not have any such concrete description at our disposal: indeed there will be no analog of theorem \ref{thm:lambdapoly} and lemma \ref{lem:diago}.

\begin{ex}\label{ex:rg} If $G$ is a group (say finite), the representation ring $R(G)$ is obtained from the monoid $R(G)^+$ of isomorphism classes of representations of $G$. We can definie operations $\l^n$ on $R(G)^+$ by setting $\l^n(V)=$ the $n$-th exterior power of $V$. It is classical that $$\l^n(V\oplus W)=\sum_{p+q=n}\l^p(V)\otimes \l^q(W).$$ So if we consider the map 
$\l_t: R(G)^+ \to 1 + tR(G)[[t]]$ defined by
$\l_t(V)=1 + \l^1(V)t + \l^2(V)t^2 + \ldots$ we see that it satisfies $\l_t(x+y)=\l_t(x)\l_t(y)$. Since $1 + tR(G)[[t]]$ is a group under multiplication of power series (and not just a monoid), we conclude that the map $\l_t$ extends to all of $R(G)$. In this way we see that $R(G)$ is a pre-$\l$-ring.
\end{ex}

\subsection{Operator rings.}

This definition will encapsulate the properties of $R^R$, as a ring with an extra operation of composition $\circ$. An {\em operator ring} is a ring $\Sigma$ with a composition $\star: \Sigma\times \Sigma\to\Sigma$ and a distinguished element $e$ ($\ne 1$) satisfying:
$$ \sigma_1\star (\sigma_2 \star \sigma_3)= (\sigma_1\star \sigma_2) \star \sigma_3$$
$$(\sigma_1 + \sigma_2)\star\sigma_3=(\sigma_1\star\sigma_3) + (\sigma_2\star\sigma_3)$$
$$(\sigma_1 \sigma_2)\star\sigma_3=(\sigma_1\star\sigma_3)  (\sigma_2\star\sigma_3)$$
$$\sigma\star e = e\star \sigma=\sigma.$$

In fact, if $\Sigma$ is an operator ring, the (adjoint) map $i: \Sigma\to \Sigma^{\Sigma}$ defined by $i(\sigma)(\tau)=\sigma\star\tau$ is a map of operator rings (in the obvious sense) which is injective. Indeed, an alternative definition of operator ring would require the existence of an injective map such as $i$ whose image is closed under the usual composition on $\Sigma^{\Sigma}$.

The following result is fundamental.

\begin{thm} The ring $\L$ is an operator ring.
\end{thm}

The reader will find in \cite{knut} the classical argument which from our point of view is rather roundabout: it is proved that $\L$ is itself a pre-$\l$-ring with good properties\footnote{in the literature one expresses these good properties by saying that $\L$ is a $\l$-ring. In our presentation this is premature, as we have not defined $\l$-rings just yet.}, which translates into a map $\L\to \L^{\L}$ as required.

A proof which is closer in spirit to our presentation can be sketched as follows. As we have observed in example \ref{ex:rg}, the elements of $\L$ act on $R(G)$, the representation ring of any finite group. Since these actions are natural in $G$, we obtain a map $\L\to F(R,R)$, where $F(R,R)$ is the ring of all natural transformations of functors between $R$ and itself. The latter is clearly an operator ring, with $\star$ being the usual composition of natural transformations. One then proves that this map is injective, and that the image of $\L$ is closed under composition. This defines an operation $\star$ on $\L$.

For the convenience of the reader, we shall indicate here how to compute $U\star V$ where $U$ and $V$ are representations of $S_m$ and $S_n$ respectively. Note that the operation $\star$ is not linear in its second argument, so we are not giving a complete description of the composition in $\L$, only partial information for the sake of concreteness. 

So let $V$ be a representation of $S_n$, and consider the endofunctor $F_V$ obtained on the category of finite dimensional, complex vector spaces by setting $$F_V(W)=(V\otimes W^{\otimes n})^{S_n}.$$ Taking $T$ to be the diagonal matrix $(t_1,\ldots,t_k)$ acting on $\C^k$, we obtain the polynomial $trace(F_V(T))$, which is homogeneous of degree $n$. This polynomial is invariant under permutation of the $t_i$'s, if only because the definition of $F_V$ is coordinate-free (the group $S_n$ plays no role here). So we obtain a polynomial $P_V$ such that $P_V(\epsilon_1,\ldots,\epsilon_k)=trace(F_V(T))$ where $\epsilon_i$ is the $i$-th elementary symmetric function in the variables $t_1,\ldots,t_k$. It is easy to see that $P_V$ is independent of $k$ as long as $k\ge n$.

To take a basic example, the functor $F_{\l^n}$ acts by $F_{\l^n}(W)=\l^n(W)$, the $n$-th exterior power of $W$ in the usual sense, and $P_{\l^n}(\epsilon_1,\epsilon_2,\ldots)=\epsilon_n$. In passing, one observes that $V\mapsto P_V$ is additive and multiplicative. Therefore, if we take an element in $\L$ which can be written $V=P(\l^1,\ldots,\l^i,\ldots)$ for a certain polynomial $P$ {\em with nonnegative coefficients}, we obtain immediately $P_V=P$. Thus $V\mapsto P_V$ (and so also $V\mapsto F_V$) is injective on the set of such elements $V$.

Motivated by this, it is then clear how to proceed. We simply observe that the argument above may be repeated with the composite functor $F_U\circ F_V$, thus yielding a polynomial $P_{U,V}$. The composite $U\star V$ is simply $P_{U,V}(\l^1,\ldots,\l^i,\ldots)$.

\begin{rmk} The proof of theorem \ref{thm:lambdapoly} is little more than this and a counting argument to show that the subring generated by the operations $\l^i$ is all of $\L$.
\end{rmk}

\subsection{A word on $\spec\L$.}

Let $X=\spec\L$ be the affine scheme associated to $\L$. The points of $X$ in the ring $R$, that is the set $X(R)=Hom_{rings}(\L,R)$, can be naturally identified with $X(R)=1 + tR[[t]]$. In other words, the functor of points of $X$, from rings to sets, can be seen as $R\mapsto 1 + tR[[t]]$.

Now, the set $1 + tR[[t]]$ can be made into an abelian group via the multiplication of power series. Since this is natural in $R$, we conclude that there is a morphism of schemes $X\times X \to X$. Looking at the coordinate rings yields a homomorphism of rings $\L\to \L\otimes\L$ which is none other than $\Delta$. So $\Delta$ is indeed a homomorphism of rings.

However, we can go further and put a ring structure on $1 + tR[[t]]$. The product of $1 + \sum a_nt^n$ and $1 + \sum b_nt^n$ is $$1 + \sum P_n(a_1,\ldots,a_n;b_1,\ldots,b_n)t^n$$ where $P_n$ is defined as follows. Consider the product $\prod (1 + X_iY_ju)$. Then $P_n$ is characterized by the fact that the coefficient of $u^n$ in this expression is $P_n(\epsilon_1,\ldots,\epsilon_n;\epsilon_1',\ldots,\epsilon_n')$, where $\epsilon_i$, resp. $\epsilon_i'$, is the $i$-th symmetric function in the $X_j$'s, resp. the $Y_j$'s. In fact, one can show that $1 + tR[[t]]$ is the ring of big Witt vectors of $R$.

It follows that $X$ is a {\em ring scheme}. In other words, $\L$ possesses another diagonal $\Delta^m:\L \to \L\otimes\L$, such that $\L$ with its two diagonals satisfies all the axioms for a commutative ring (associativity, distributivity, etc), only with all the arrows reversed.

\subsection{The definition of a $\l$-ring.}

\begin{defn} A {\em $\l$-ring} is a pre-$\l$-ring $R$ such that:
\begin{enumerate}
\item the map $\theta: \L\to R^R$ is a map of operator rings,
\item the following diagram commutes:
$$\begin{CD}
    \L @>{\Delta^m}>> \L\otimes\L
\\ @V{\theta}VV        @VV{\theta^2}V
\\  R^R @>m>> R^{R\times R} 
\end{CD}$$ Here $m(f)(x,y)=f(x\cdot y)$ ($m$ is for multiplication).
\end{enumerate}
\end{defn}
From the second axiom, one can derive a formula for $\l^n(x\cdot y)$, just as we did in the case of pre-$\l$-rings for $\l^n(x+y)$: $$\l^n(x\cdot y)=P_n(\l^1(x),\ldots,\l^n(x);\l^1(y),\ldots,\l^n(y)).$$ Again this formula is equivalent to the commutativity of the diagram. Similarly, we could replace the first axiom by a single formula involving universal polynomials. Details in any of \cite{knut}, \cite{atiyahtall}, \cite{patras}.

\begin{ex} The ring $R(G)$ as in example \ref{ex:rg} is a $\l$-ring. In fact, considering the definition of the operation $\star$ on $\L$, the composition axiom for $R(G)$ is tautological. As for the multiplicative axiom, one way around it is to study the restrictions to abelian subgroups of $G$, for which any representation is a direct sum of $1$-dimensional representations; for such an element the operations $\l^n$ are $0$ as soon as $n\ge 2$, and then the required formulae are easily established. See \cite{atiyahtall} for the most concise proof.
\end{ex}

\begin{ex}{\em The free $\l$-ring on one generator.} The ring $\L$ itself may be viewed as a $\l$-ring. The map $\L\to \L^{\L}$  is the adjoint map sending $a$ to the map $b\mapsto a\star b$. As already pointed out, this is a map of operator rings. As for the additive and multiplicative axioms, we only need to notice that they hold in $R(G)$ for any $G$, and then embed $\L$ in $F(R,R)$ as above to conclude that they hold in $\L$.

With this structure, we have $\l^n(\l^1)=\l^n$ (that is $\l^n\star\l^1=\l^n$, which is plainly true). Since there are no algebraic relations between the $\l^n$'s, $\L$ is called the free $\l$-ring on one generator, that generator being $\l^1$. Another justification is that if we pick $x$ in any $\l$-ring $R$, then there is a unique map $\L\to R$ which sends $\l^1$ to $x$ and commutes with the operations. 
\end{ex}

\section{$\b$-rings}\label{section:greek2}

\subsection{The ring $\B$, pre-$\b$-rings and $\b$-rings.}

We let $$\B=\bigoplus_{n\ge 0} A(S_n)$$ where here and elsewhere $A(G)$ denotes the Burnside ring of the finite group $G$. ($\B$ is meant to be a capital $\b$.) We define a multiplication on $\B$ by setting $$\rho\sigma=Ind_{S_p\times S_q}^{S_{p+q}} (\rho\times\sigma).$$ Thus $\B$ becomes a commutative ring with unit.

We further define $$\B^2=\bigoplus_{p, q} A(S_p\times S_q).$$ 

There is a natural map $\Delta: \B \to \B^2$ defined using the restriction maps $A(S_n) \to A(S_p\times S_q)$ (for $p+q=n$). 

For any finite group $G$, the Burside ring $A(G)$ is a free abelian group, a canonical basis being given by the isomorphism classes of transitive $G$-sets, or equivalently the conjugacy classes of subgroups of $G$. Thus if $H$ is a subgroup of $S_n$, resp $S_p\times S_q$, we get an element $\b_H\in \B$, resp $\b_H\in\B^2$, corresponding to the conjugacy class of $H$. In this notation, the multiplication is given by $\b_H\b_K=\b_{H\times K}$.

\begin{defn}\label{defn:bring} A {\em pre-$\b$-ring} is a ring $R$ together with two maps $\theta: \B\to R^R$ and $\theta^2: \B^2\to R^{R\times R}$ such that 
\begin{enumerate}
\item if $A$, resp $B$, is a subgroup of $S_p$, resp $S_q$, we have
$$\theta^2(\b_{A\times B})(x,y)=\theta(\b_A)(x)\cdot \theta(\b_B)(y),$$ 

\item the following diagram commutes:
$$\begin{CD}
    \B @>{\Delta}>> \B^2
\\ @V{\theta}VV        @VV{\theta^2}V
\\  R^R @>a>> R^{R\times R} 
\end{CD}$$
\end{enumerate}
\end{defn}

We shall often suppress $\theta$ and $\theta^2$ from the notations, thus viewing $\b_H\in \B$ as a map $\b_H: R\to R$, and similarly for elements of $\B^2$. In this way the first condition above is simply $\b_{A\times B}(x,y)=\b_A(x)\cdot\b_B(y)$.

\begin{ex}\label{ex:ag} If $G$ is a group (say finite), the Burnside ring $A(G)$ is obtained from the monoid $A(G)^+$ of isomorphism classes of set-theoretic representations of $G$. We can definie operations $\b_H$ on $A(G)^+$, for $H$ a subgroup of $S_n$, by setting $\b_H(X)=X^n/H$.

In order to extend this to an operation on $A(G)$, the quickest way is to follow Rymer \cite{rymer} and use tom Dieck's definition of the Burnside ring. Suppose we define an equivalence relation on the set of all reasonable $G$-spaces, by identifying $X$ and $Y$ if for any subgroup $H$ of $G$, the spaces of fixed points $X^H$ and $Y^H$ have the same Euler characteristic. Then the set of equivalence classes is a ring under disjoint union and cartesian product, which is isomorphic to $A(G)$ (details can be found in \cite{tammo}). The expression $X^n/H$ makes sense for any $G$-space $X$, so we have indeed an operation $\b_H$ on $A(G)$.

In order to define now the map $\theta^2$, we need to specify for each subgroup $H$ of $S_p\times S_q$ an operation $\b_H: A(G)\times A(G)\to A(G)$. We set:
$$\b_H(X,Y)=(X^p\times Y^q)/H$$ for the $G$-spaces $X$ and $Y$. 

We are requested to prove the addition axiom, that is the equality
$$\b_H(x+y)=\Delta(\b_H)(x,y)$$ for $H$ a subgroup of $S_n$ and $x,y\in A(G)$. All the other examples of the additivity axiom in this article will follow from this one, so we present the details.

\begin{proof} For $G$-spaces $X$ and $Y$, we start by writing 
$$\frac{(X\coprod Y)^n}{H}=\coprod_{p+q=n} A^{p,q}$$ where $A^{p,q}$ is the set of classes of $n$-tuples in $(X\coprod Y)^n$ for which exactly $p$ points are taken from $X$. We claim that $$A^{p,q}=(S_n/H)_{|S_p\times S_q}(X,Y).$$ To makes sense of this we recall that the notation $(S_n/H)_{|S_p\times S_q}$ refers to the restriction of the $S_n$-set $S_n/H$ to the subgroup $S_p\times S_q$. Moreover this is seen as an element of $\B^2$, and the map $\theta^2$ is suppressed. Clearly the claim is all we need to establish.

The orbits of $S_p\times S_q$ in $S_n/H$ correspond to the double cosets $(S_p\times S_q)\sigma H$. Let us pick one. The stabilizer of $\sigma H$ is $(S_p\times S_q)\cap \sigma H \sigma^{-1}$. We define a map 
$$\begin{CD}
\iota_\sigma: \frac{(X^p\times Y^q)}{(S_p\times S_q)\cap \sigma H \sigma^{-1}} @>>> A^{p,q} \end{CD}$$ by $$(x_1,\ldots,x_p;y_1,\ldots,y_q)\mapsto (x_1,\ldots,x_p;y_1,\ldots,y_q)\cdot \sigma^{-1}$$ where the action of $\sigma^{-1}$ on the right is, as usual, by permuting the entries of the $n$-tuple. This is well-defined, and it is not difficult to see that $\iota_\sigma$ is injective (essentially because the union of $X$ and $Y$ is disjoint).

Further, if $\iota_\sigma$ and $\iota_\tau$ do not have disjoint images, then we see that $$\iota_\sigma(\{1,2,\ldots,p\})=\iota_\tau(\{1,2,\ldots,p\})$$ so that $\tau\sigma^{-1}\in S_p\times S_q$. It follows that we have an injective map 
$$\begin{CD}
\coprod_\sigma \frac{(X^p\times Y^q)}{(S_p\times S_q)\cap \sigma H \sigma^{-1}} @>{\coprod \iota_\sigma}>> A^{p,q}
\end{CD}$$ if we only take $\sigma$ to run over representatives for the double cosets. The latter is plainly surjective. This is the desired isomorphism of $G$-spaces.\end{proof}

\end{ex}

The following result is fundamental.

\begin{thm}\label{thm:bopring} The ring $\B$ is an operator ring.
\end{thm}

A proof of this theorem has been given by Vallejo \cite{vallejofree}, and it may be sketched as follows. As we have observed in example \ref{ex:ag}, the elements of $\B$ act on $A(G)$, the Burnside ring of any finite group. Since these actions are natural in $G$, we obtain a map $\B\to F(A,A)$, where $F(A,A)$ is the ring of all natural transformations of functors between $A$ and itself. The latter is clearly an operator ring, with $\star$ being the usual composition $\circ$ of natural transformations. One then proves that this map is injective, and that the image of $\B$ is closed under composition. This defines an operation $\star$ on $\B$.

The operation $\star$ satisfies in particular $\b_H\star\b_K=\b_{H\wr K}$. Since $\star$ is not linear in its second variable, this does not describe the composition completely.

\begin{defn} A {\em $\b$-ring} is a pre-$\b$-ring such that the map $\theta: \B\to R^R$ is a map of operator rings.
\end{defn}

There is no other requirement for $\b$-rings: see the comments at the end of this section.

\begin{ex} The ring $A(G)$ as in example \ref{ex:ag} is a $\b$-ring. In fact, considering the definition of the operation $\star$ on $\B$, the composition axiom for $A(G)$ is tautological.
\end{ex}

\begin{ex}{\em The free $\b$-ring on one generator.} The ring $\B$ itself may be viewed as a $\b$-ring. The map $\B\to \B^{\B}$  is the adjoint map sending $a$ to the map $b\mapsto a\star b$. As already pointed out, this is a map of operator rings. 

To define the map $\theta^2$, we will need to see elements of $\B$ as natural transformations in $F(A,A)$ as above. Given an element $\b_H\in\B^2$, we have then a natural transformation $\b_H: A(G)\times A(G)\to A(G)$ as in example \ref{ex:ag}. If we take $a$ and $b$ in $\B$, we can consider the transformation $A(G)\to A(G)$ given by $x\mapsto \b_H(a(x),b(x))$. It is not hard to see that this transformation is in $\B$ (in fact this statement follows from lemma \ref{lem:bplus} below) and we take it for $\b_H(a,b)$, which as usual is short hand for $\theta^2(\b_H)(a,b)$.

Since the additive axiom holds for $A(G)$ for any $G$, it follows that it also holds for $\B$.

With this structure, we have $\b_H(\b^1)=\b_H$ (that is $\b_H\star\b^1=\b_H$, which is plainly true). Therefore, $\B$ is generated by $\b^1$ as a $\b$-ring, and it is in fact called the free $\l$-ring on one generator. A justification is that if we pick $x$ in any $\b$-ring $R$, then there is a unique map $\B\to R$ which sends $\b^1$ to $x$ and commutes with the operations. 
\end{ex}

\subsection{Polynomial operations.}

It will be useful in the sequel to have an easy criterion to check whether a given pre-$\b$-ring is a $\b$-ring or not. It will prove convenient to use {\em polynomial operations}, also known as maps of finite degree, or algebraic maps. A map $f: A \to B$ between the abelian group or monoid $A$ and the abelian group $B$ is said to have degree $\le 0$ if it is constant. Inductively, $f$ is said to have degree $\le n$ if for all $a\in A$, the map $x\mapsto f(x + a) - f(x) -f(a)$ has degree $\le n-1$. The degree of $f$ is the smallest $n$ such that $f$ has degree $\le n$, if it exists at all in which case we say that $f$ has finite degree (or is polynomial).

Polynomial maps have the following useful properties.

\begin{lem}\label{lem:polyagree} Suppose that $A$ is the group completion of the monoid $A^+$. Then any polynomial map $A^+\to B$ extends uniquely to a polynomial map $A\to B$.

In particular, if two polynomial maps $A\to B$ agree on a family of generators of $A$ which is closed under addition, then they are equal.
\end{lem}

This is easily proved by induction on the degree.

The relevance to $\b$-rings is the following:

\begin{lem}\label{lem:composepoly} For any $a\in\B$, the map $b\mapsto a\star b$ is polynomial.
\end{lem}

\begin{proof} The point is that the operations $\b_H$ on Burnside rings as in example \ref{ex:ag} are polynomial maps, as was proved by Vallejo \cite{vallejopoly}. If we inject $\B$ into $F(A,A)$ as explained after theorem \ref{thm:bopring}, the result becomes obvious. 
\end{proof}

Let us introduce the additive monoid $\B^+\subset \B$ consisting of all elements of the form $\sum n_K\b_K$ with each $n_K$ nonnegative (that is, $\B^+$ consists of all sums of $S_n$-sets, for various $n$). Similarly, we have the monoid $(\B^2)^+$, and the diagonal $\Delta$ maps $\B^+$ into $(\B^2)^+$. Then, as an application of the above, we have:

\begin{coro}\label{coro:bplus} Let $R$ be a pre-$\b$-ring, and suppose that all the operations $\b_H:R\to R$ are polynomial. Then $R$ is a $\b$-ring if and only if for any subgroup $H$ of $S_n$ and any element $b\in\B^+$, we have
$$(\b_H\star b)(x)=\b_H(b(x)),$$ for all $x\in R$.
\end{coro}

\begin{proof} From the definitions, $R$ is a $\b$-ring if and only if for all $x\in R$, $a, b\in \B$, the relation $$(a\star b)(x)=a(b(x))$$ holds. Since $\star$ is linear in its first variable, it suffices to check this for $a=\b_H$ for all $H$. Now, both sides of the equation in the statement of the corollary are polynomial maps of $b$, so from the lemma it suffices to check that they are equal on $\B^+$.
\end{proof}

Strictly speaking, we shall not use this lemma in this form (though in the course of the proof of theorem \ref{thm:pizero} we shall use an infinitesimal variant). It is included to illustrate the simplification that polynomial operations may provide.

Finally we record:

\begin{lem}\label{lem:bplus} $\B^+$ is closed under $\star$.
\end{lem}

\begin{proof} This follows from \cite{vallejofree}, proposition 2.1.
\end{proof}

\subsection{Relations between $\B$ and $\L$.}

For any finite group $G$, there is a natural map $A(G)\to R(G)$ obtained by associating to a $G$-set $X$ the complex vector space with basis $X$. Thus we have maps $\B\to\L$ and $\B^2\to \L\otimes\L$. The following is obvious.

\begin{prop} The following diagram commutes.
$$\begin{CD}
	\B @>{\Delta}>> \B^2
\\  @VVV            @VVV
\\  \L @>{\Delta}>> \L\otimes\L
\end{CD}$$

Moreover, the natural map $\B\to\L$ is a map of operator rings. Also, it is surjective.
\end{prop}

If $R$ is a pre-$\l$-ring with associated map $\theta: \L\to R^R$, we may consider the composition
$$\begin{CD} 
\theta': \B @>>> \L @>{\theta}>> R^R
\end{CD}~.$$
Thanks to the above theorem, we see that $R$ becomes in this way a pre-$\b$-ring, which we shall denote by $\b(R)$. Moreover, if $R$ is a $\l$-ring, then $\b(R)$ is a $\b$-ring. (Conversely, it is not hard to see that when $R$ is not a $\l$-ring because the composition axiom fails, then $\b(R)$ is not a $\b$-ring.)

We may also proceed the other way around. Consider the element $\b_{S_n}\in\B$. It corresponds to the trivial $S_n$-set, so its image under $\B\to\L$ is $\b^n$. We pause to indicate some notations. 

\begin{notations2} In the sequel, we shall write $\b^n$ to denote the element $\b_{S_n}$ in $\B$, as well as its image in $\L$. Moreover, we shall also write $\b_H$ for the element in $\L$ which is the image of the previously defined element $\b_H\in\B$. No confusion should arise.
\end{notations2}

Given a pre-$\b$-ring $R$, we have operations $\b^n: R\to R$ satisfying $$\b^n(x+y)=\sum_{p+q=n}\b^p(x)\cdot \b^q(y).$$ Therefore, we may give $R$ the structure of a pre-$\l$-ring which we denote by $\l(R)$. However, we caution that $\l(R)$ may very well fail to be a $\l$-ring even if $R$ is a $\b$-ring.

Put differently, we have a section $s:\L\to\B$ which behaves well with respect to the diagonals, but it is not compatible with $\star$. In fact, the subring of $\B$ generated by the $\b^n$'s is not closed under composition.

It is clear from the definitions that $\l(\b(R))$ is isomorphic to $R$ as a (pre-)$\l$-ring. However, $\b(\l(R))$ is very different from $R$: for example, in $\b(\l(R))$, any operation in the kernel of $\B\to\L$ acts as $0$, which of course is not the case in general. So there are more $\b$-rings than $\l$-rings.

\begin{ex} The pre-$\l$-ring $\l(A(G))$ is not a $\l$-ring unless $G$ is cyclic: Siebeneicher \cite{sieben} has proved that the multiplicative axiom fails, while Gay, Morris and Morris \cite{gay} have proved that the composition axiom fails (in either case they check this on Adams operations -- see below).

This pre-$\l$-ring structure on $A(G)$ is in competition with another one, in which we take $\l^n(X)=$ the set of subsets of $X$ of cardinal $n$. These two structures are different unless the order of $G$ is odd (\cite{sieben}). Moreover, the second structure we have introduced is not in general compatible with the $\l$-structure on $R(G)$, in that the natural map $A(G)\to R(G)$ is not in general a map of pre-$\l$-rings (consider $G=\Z/2$ for example). However, this is the case with $\l(A(G))$.
\end{ex}

Let us say a word about the way in which $\b$-rings arise in practice. One often starts with a ring $R$ which possesses natural operations $\b^n$ (or $\l^n$), satisfying the right formula for $\b^n(x+y)$ (or $\l^n(x+y)$). Taking this as a pre-$\l$-ring structure on $R$, one checks if we have a $\l$-ring or not. If not, the second best thing one can do is consider the operations $\b^n$ as operations $\b_{S_n}$ and try to find the correct definition for operations $\b_H$, for $H$ any subgroup of $S_n$. If one succeeds, with a bit of luck one can define a $\b$-ring structure on $R$, which at least affords a way of computing the composition of operations in $R$.

\subsection{Comments.}
\begin{enumerate}
\item The most glaring difference between $\l$-rings and $\b$-rings as we define them is the lack of a compatibility requirement between the operations $\b_H$ and the products in the ring -- that is, we do not insist on a formula for $\b_H(x\cdot y)$. In fact, the current definition of $\b$-rings in the literature does not even include a condition on the behavior of the operations with respect to addition. We propose to impose the commutativity of the diagram above, and the very existence of the map $\theta^2$, because it seems a natural analog of the formulae which one can write down in the case of $\l$-rings. This condition is satisfied by all the known examples of $\b$-rings.

The difficulty lies, of course, in the understanding of the ring $\B$. For example, I am not aware of any known minimal set of multiplicative generators. I do not know how to construct an analog of the second diagonal $\Delta^m$.

\item Let us deliver a word of caution, for a number of mistakes have appeared in the literature. First, it is important to keep in mind that the $\b$-ring structure on $\B$  is {\em not} the direct sum of the structures of the various $\b$-rings $A(S_n)$. Indeed, the $\b_H$ operations on $\B$ do not preserve degrees. Likewise for the $\l$-ring structure on $\L$. I am aware of two unfortunate printed occurences of this confusion.

Similarly, a couple of enthusiastic authors have used a version of corollary \ref{coro:bplus} in which $b$ is only taken to be of the form $\b_K$; as such the result is not true, and the confusion comes from forgetting the words "closed under addition" from lemma \ref{lem:polyagree}.

\item Here is a possible variant of the definition of $\b$-ring. One may want to define $\B^n$ using $n$-fold products of symmetric groups, and to require the existence of maps $\B^n\to\B^m$ when $n\le m$ (possibly organised in some standard algebraic structure), giving formulae for $\b_H(x_1 + x_2 + \ldots x_n)$. In a given degree, $\B^n$ becomes trivial as $n$ grows large, and it would follow automatically that the operations $\b_H$ on any $\b$-ring are polynomial. In practice this is quite useful to know.
\end{enumerate}

\section{$\pi^0(X)$ is a $\b$-ring}

\subsection{The theorem of Priddy-Quillen-Segal.}

A celebrated theorem, due to these three authors as well as others, asserts that $QS^0$ has the homotopy-type of $\zbs$, so that 
$$\pi^0(X)=[X,\zbs]$$ for any space $X$. It is therefore expected that $\pi^0(X)$ should be related to the monoid $Cov(X)^+$ of (isomorphism classes of) covering spaces of $X$, much as $K^0(X)$ is related to the complex vector bundles over $X$. We have the following precise formulation, using the natural transformation of functors $Cov^+\to \pi^0$.

\begin{thm}\label{thm:segal} Let $T$ be any abelian-group valued, representable, homotopy functor. Then any polynomial transformation of functors $Cov^+\to T$ extends to a unique (polynomial) transformation $\pi^0 \to T$.
\end{thm}

By polynomial transformation we mean that each map $Cov(X)^+\to T(X)$ is polynomial. This theorem was first stated by Segal with "additive" replacing "polynomial" \cite{segal}, and Vallejo proved the version which we state here \cite{vallejopizero} (essentially by using lemma \ref{lem:polyagree} at appropriate places).

At this point it may be useful to recall how the transformation $Cov(X)^+\to \pi^0(X)$ is obtained. Assume that $X$ is connected -- if not, proceed by component. Given an $n$-sheeted covering $E\to X$, there is a uniquely defined homotopy class of maps $f:X\to BS_n$ such that $E=f^*(E_n)$ (the pullback covering), where $E_n\to BS_n$ is the "universal" $n$-sheeted cover of $S_n$. The latter is not to be mistaken with $ES_n\to BS_n$, the "universal principal cover", which has degree $n!$ (the difference between them is exactly analogous to that between the universal vector bundle over $BU(n)$ and the universal $U(n)$-principal bundle over $BU(n)$).

Composing $f$ with the natural map $BS_n\to \{n\}\times \bs\to \zbs$, one obtains a class in $\pi^0(X)$, as required.

Now let $Cov(X)$ denote the group obtained from the monoid $Cov(X)^+$. Then we have the factorisation $Cov(X)^+ \to Cov(X) \to \pi^0(X)$.

\begin{ex} When $X=BG$ is the classifying space of a finite group $G$, we have $Cov(X)^+=A(G)^+$ and $Cov(X)=A(G)$. Thus in this case $Cov(X)$ is already a $\b$-ring. We shall see that the operations on the Burnside ring extend to operations on $\pi^0(BG)$, which by the Segal conjecture, now a theorem, is the completion of $A(G)$ at the augmentation ideal.
\end{ex}

Finally we note the following easy corollary:

\begin{coro}\label{coro:segal} Under the same hypotheses on $T$, any natural transformation $$Cov^+ \times Cov^+ \to T$$ extends uniquely to a natural transformation
$$\pi^0 \times \pi^0 \to T.$$
\end{coro}

To be precise, $Cov^+ \times Cov^+$ means the functor of two variables $X,Y\mapsto Cov(X)^+\times Cov(Y)^+$, while $T$ is a functor of two variables with is assumed to satisfy the hypotheses of theorem \ref{thm:segal} when either of the variables is fixed.

\subsection{The operations.}

Let us fix a base space $X$. We shall first introduce operations $\b_H$ on $Cov(X)^+$ by setting $$\b_H(E)=E^n/H$$ when $H$ is a subgroup of $S_n$. Here $E^n$ denotes the $n$-fold fibred product of $E$ with itself, over $X$: $$E^n=E\times_X E \times_X \cdots \times_X E \qquad (n~ \textrm{factors}).$$

We see $\b_H$ as a map $Cov(X)^+\to Cov(X)$. Let us prove that it is polynomial. For this, we consider the case $X=BG$ first. Then $\b_H$ can be seen as the operation $A(G)^+ \to A(G)$ considered in \ref{ex:ag}. As indicated there, this is a polynomial operation; in fact, if $H$ is a subgroup of $S_n$, then $\b_H$ has degree $n$ (\cite{vallejopoly}, 5.6). Note that this is independent of $G$.

For the general case we shall use the following lemma.

\begin{lem}\label{lem:pullback}
Let $X$ be connected, and let $E_1, E_2, \ldots, E_k$ be elements of $Cov(X)^+$. Then there is a finite group $G$ and a single map $f: X \to BG$ such that $E_i=f^*(F_i)$ for a cover $F_i$ of $BG$, for all $i$ simultaneously.
\end{lem}

\begin{proof} Let $f_i: X \to BS_{n_i}$ be a map classifying $E_i$, ie $E_i=f_i^*(E_{n_i})$. Put $G=S_{n_1}\times S_{n_2}\times \cdots \times S_{n_k}$, let $p_i: G \to S_{n_i}$ be the projection, and let $F_i=p_i^*(E_{n_i})$. Finally, let $f=f_1\times\cdots\times f_k$, so that $p_i\circ f=f_i$.
\end{proof}

Now we shall prove that $\b_H: Cov(X)^+ \to Cov(X)$ has degree $n$ also (where we assume that $X$ is connected without loss of generality). Indeed, this means that a certain function $\widehat{\b_H}$ of $n$ variables taken in $Cov(X)^+$ is constant; morever $\widehat{\b_H}$ is defined from $\b_H$ in algebraic terms. In particular, $\widehat{\b_H}$ is natural in $X$, just as $\b_H$ is. However, any $n$ elements in $Cov(X)^+$ can be seen as the pullbacks of $n$ elements in $Cov(BG)^+$ for a certain group $G$, by the lemma. From the case just considered, we see that $\widehat{\b_H}$ is constant on $Cov(BG)^+$, and from the definitions its constant value must be $0$ (evaluate on the empty covers). Therefore by naturality, $\widehat{\b_H}$ is $0$ on $Cov(X)^+$ as well.

Composing with the map $Cov(X)\to \pi^0(X)$ which is additive, we see that $\b_H$ is also polynomial when seen as a map $Cov(X)^+\to \pi^0(X)$.

Appealing to theorem \ref{thm:segal}, we may write the following commutative diagram:
$$\begin{CD}
	Cov(X)^+ @>{\b_H}>> Cov(X)^+
\\    @VVV                @VVV
\\  \pi^0(X) @>{\b_H}>> \pi^0(X) 
\end{CD}$$
In practice, naming both horizontal arrows $\b_H$ will cause no confusion.

Similarly, when $H$ is a subgroup of $S_p\times S_q$, we define an operation $\b_H: Cov(X)^+ \times Cov(Y)^+ \to Cov(X\times Y)^+$ by setting $$\b_H(E,F)=(E^p\times F^q)/H.$$ In this expression, we see $E$ resp $F$ as a cover of $X\times Y$ by pulling back along the projection onto $X$ resp $Y$. Then the products are fibred products as above. 

Composing with $Cov(X\times Y)^+ \to \pi^0(X\times Y)$ and using corollary \ref{coro:segal}, whose hypotheses are satisfied as one sees easily, we end up with a natural transformation $\pi^0(X)\times \pi^0(Y) \to \pi^0(X\times Y)$. There is also a commutative diagram as before. When $X=Y$, we may use the diagonal $X\to X\times X$ to pull back from $\pi^0(X\times X)$ to $\pi^0(X)$. The bottom line is that we have a natural transformation $$\b_H: \pi^0(X)\times \pi^0(X) \to \pi^0(X).$$ We may also use the notation $\b_H$ for the transformation on the level of $Cov^+$, that is $\b_H(E,F)=(E^p\times_X F^q)/H\in Cov(X)^+$ when $E$ and $F$ are in $Cov(X)^+$.

Thus we have defined two maps $\theta: \B\to \pi^0(X)^{\pi^0(X)}$ and $\theta^2: \B^2 \to \pi^0(X)^{\pi^0(X)\times \pi^0(X)}$ sending $\b_H$ to the operation defined above. In what follows we shall prove that $\pi^0(X)$ endowed with these is a $\b$-ring.

\subsection{The axioms.}

We may now prove:

\begin{thm}\label{thm:pizero} $\pi^0(X)$ is a $\b$-ring.
\end{thm}

\begin{proof}

We start with the composition axiom; that is let us prove, for $x\in\pi^0(X)$ and $a,b\in\B$, that $$(a\star b)(x)=a(b(x)).$$ Now either side of this equation defines a polynomial transformation $Cov(X)^+ \to \pi^0(X)$ and by theorem \ref{thm:segal} it is sufficient to show that these agree. In other words, we may restrict our attention to $x\in Cov(X)^+$. We reason then as in the proof of corollary \ref{coro:bplus}: we may reduce the problem to $a=\b_H$ for all $H$ by linearity, and the resulting equation to be proved is a polynomial function of $b$ (using lemma \ref{lem:composepoly}), so we may restrict attention to $b\in\B^+$.

When $c\in\B^+$ and $x\in Cov(X)^+$, we may see $c(x)$ as an element of $Cov(X)^+$. Applying this with $c=a=\b_H$, $c=b$, and then $c=a\star b$ (using \ref{lem:bplus}), we conclude that the equation $(\b_H\star b)(x)=\b_H(b(x))$ makes sense in $Cov(X)^+$, and if we can prove it to be true, then we are done. Indeed, both sides of the equation are natural in $X$; the equation is true when $X=BG$ since $A(G)$ is a $\b$-ring; so we may use the same trick again.

It remains to prove the addition axiom. Recall from definition \ref{defn:bring} that we are required to prove
$$\b_H(x+y)=\Delta(\b_H)(x,y)$$ for all $H$ and all $x,y\in\pi^0(X)$.

Let $x$ resp $y$ be a cover of $X$ resp $Y$. We may see $x$ and $y$ as covers of $X\times Y$ by pulling back along the projections. Recall that for any element $b\in(\B^2)^+$ we have defined $b(x,y)\in Cov(X\times Y)^+$. Thus it makes sense to ask for the equation $\b_H(x+y)=\Delta(\b_H)(x,y)$ to hold in $Cov(X \times Y)^+$. Indeed, if $X$ and $Y$ are replaced by $BG_1$ and $BG_2$ respectively, the formula does hold since $A(G_1\times G_2)$ is a $\b$-ring. By naturality again, it follows that the formula holds in general.

Using corollary \ref{coro:segal}, we obtain an equality of natural transformations $$\pi^0(X)\times \pi^0(Y) \to \pi^0(X\times Y).$$ Taking $X=Y$ and pulling back along the diagonal gives the result.
\end{proof}

\begin{ex} When $X$ is reduced to a point, we have $\pi^0(*)=Cov(*)=\Z$. Now, $\Z$ has a unique structure of $\l$-ring in which $\l^k(r)=\left(  {r \atop k}  \right)$. The $\b$-ring $\pi^0(*)$ is $\b(\Z)$ (notations as in section \ref{section:greek2}). Explicitly we have, for $H$ a subgroup of $S_n$:
$$ \b_H(r)=\frac{1}{| H |} \sum_{k=1}^n c_k r^k$$ where $c_k$ is the number of elements in $H$ which consist of a product of $k$ disjoint cycles.
\end{ex}

\begin{rmk}\label{rmk:reduced} The action of $\b_H$ is natural in $X$, and $\b_H(0)=0$ in $\pi^0(*)$, so the augmentation ideal
$$\tilde{\pi}^0(X)=\ker ( \pi^0(X) \to \pi^0(*) )$$ is stable under the action of $\B$. We call such an object (that is, a $\b$-ring without unit) a $\b$-ideal.

For example, when $X=S^n$, then $\tilde{\pi}^0(S^n)=\pi_n^s(S^0)$, the $n$-th stable homotopy group of $S^0$. Thus the stable homotopy groups of spheres are $\b$-ideals.

Let us point out that the functor $\tilde{\pi}^0$ is represented in the pointed category by $\zbs$, that is
$$\tilde{\pi}^0(X)=[X,\zbs]^\bullet ,$$ the set of pointed maps up to pointed homotopy. On the category of connected, pointed spaces it is represented by $\bs$. As $\tilde{\pi}^0(X)$ is a $\b$-ideal, we have pointed maps $\b_H: \bs \to \bs$. Much of the rest of this paper is devoted to studying these, in particular their effect on the homotopy and cohomology of $\bs$.

\end{rmk}

\subsection{Comparison with algebraic $K$-theory.}

Theorem \ref{thm:segal} is a statement about the group completion $\coprod BS_n \to \zbs$, and it holds in greater generality (see \cite{vallejopizero}, 1.10 again). Another familiar example is the group completion
$$\coprod BGL_n(A) \to \kba$$ where $A$ is any commutative ring. Here the left hand side is the (realization of the) nerve of the category of free, finitely generated $A$-modules, while the right hand side represents the algebraic $K$-theory over $A$ of spaces, that is $$K_A(X)=[X,\kba].$$ Correspondingly, the reduced algebraic $K$-theory of the connected, pointed space $X$ is
$$\tilde{K}_A(X)=[X,\ba]^\bullet.$$ 
The ($0$-th) cohomotopy of $X$ can be thought of as the algebraic $K$-theory of $X$ over "the field with one element".

The work we have done for $\pi^0$ may be repeated with $K_A$. Let us assume that $A$ is a field to simplify the discussion. One needs to replace covering spaces with (unpointed) maps $X\to \coprod BGL_n(A)$, and when $X=BG$ for a discrete group $G$ this is the monoid of isomorphism classes of finitely generated $G$-modules. The group obtained from this monoid is $R_A(G)$, the representation ring over $A$, and it is always a $\l$-ring. The rest of the argument is completely similar to the above, and we get:

\begin{thm}
When $A$ is a field, $K_A(X)$ is a $\l$-ring, and $\tilde{K}_A(X)$ is a $\l$-ideal.
\end{thm}

That $\tilde{K}_A(X)$ is a $\l$-ideal (even when $A$ is only a commutative ring) was proved in \cite{kra} using a more explicit method. In fact, the two structures of $\l$-ideal must coincide, since they do for classifying space $BG$ as is easily checked.

There is also a natural transformation $\pi^0(X) \to K_A(X)$. One way to think about it is to consider the natural map $\zbs\to \kba$ which is induced by the inclusion $S_\infty\to GL_\infty(A)$ and the map $\Z\to K_0(A)$ sending $d$ to $A^d$ (the class of the free $A$-module of rank $d$). Alternatively, we may use theorem \ref{thm:segal} again. For any $n$-sheeted cover $E$ of $X$, consider a classifying map $X\to BS_n$ and compose with $BS_n \to BGL_n(A) \to \{n\} \times \ba$. This defines a transformation $Cov(X)^+ \to K_A(X)$ which extends to our transformation $\pi^0(X) \to K_A(X)$.


We obtain easily:
\begin{prop}\label{prop:ktheory} For any commutative ring $A$ and any connected space $X$, the natural map $$\tilde{\pi}^0(X) \to \tilde{K}_A(X)$$ is a map of $\b$-ideals.
\end{prop}

\section{$\P$-rings}

\subsection{Pre-$\P$-rings and $\P$-rings.}

A {\em pre-$\P$-ring} is a ring $R$ together with additive operations $\P^k: R\to R$, for $k\ge 0$. A {\em $\P$-ring} is a pre-$\P$-ring $R$ such that the $\P$-operations are multiplicative and satisfy $$\P^k\circ\P^l=\P^{kl}.$$ (In particular, in this case the operations commute.) The operations $\P^k$ are usually referred to as the Adams operations.

\begin{ex}\label{ex:psigraded} Suppose $H^*$ is a graded ring. Then $H^*$ becomes a $\P$-ring when we set $$\P^k(x)=k^{|x|}x$$ for homogeneous $x$, and extend by linearity.
\end{ex}

\subsection{$\P$-operations in a $\l$-ring.}

 As before, we let
$$\l_t=1 + \l^1t + \l^2t^2 + \ldots\in 1 + t\L[[t]]$$
and
$$\b_t=1 + \b^1t + \b^2t^2 + \ldots\in 1 + t\L[[t]]$$ so that $\l_t\b_{-t}=1$. We define $$\P_t=-t\frac{\l_{-t}'}{\l_{-t}}=t\frac{\b_{t}'}{\b_{t}}.$$
The coefficients of $\P_t=\P^1t + \P^2t^2 + \ldots$ define elements $\P^k\in\L$.

Now, let $R$ be a pre-$\l$-ring. The relation $\l_t(x+y)=\l_t(x)\l_t(y)$ yields immediately $\P_t(x+y)=\P_t(x)+\P_t(y)$. In other words, {\em every pre-$\l$-ring is also a pre-$\P$-ring}.

Similarly, one verifies that if $R$ is a $\l$-ring, then it becomes also a $\P$-ring under the action of the operation $\P^k\in\L$ that we have just defined. Details in \cite{knut}.

\begin{ex}\label{ex:adamsfinfield} Let $\F_q$ denote the field with $q$ elements. The algebraic $K$-theory $\tilde{K}_{\F_q}(X)$ of the space $X$, as mentioned in the previous section, is a $\l$-ideal. Thus there are maps $\P^k: BGL_\infty(\F_q)^+ \to BGL_\infty(\F_q)^+$ representing the Adams operations. The effect of these maps on the homotopy groups of $BGL_\infty(\F_q)^+$, or equivalently on the algebraic $K$-theory of spheres, is known \cite{kra}: one has $\pi_{2i}(BGL_\infty(\F_q)^+)=0$ and on $\pi_{2i - 1}(BGL_\infty(\F_q)^+)=\Z/(q^i - 1)$, the map $\P^k_*$ is multiplication by $k^i$.

\end{ex}

\subsection{The ring $\L\otimes \Q$.}

If one decides to write down the equations defining the various $\P^k$ in terms of, say, the elements $\b^i$, one obtains a system of equations which can be solved {\em at the price of tensoring with $\Q$}. In other words, each $\b^i$ can be written as a linear combination with rational coefficients of the elements $\P^k$. So let us have a closer look at the ring $\L\otimes\Q$.

The rank of $\L$ in degree $n$, ie the rank of $R(S_n)$, is the number of conjugacy classes in $S_n$, and in turn this equals the number of {\em partitions} of $n$: a partition for us will be a sequence $n_1\ge n_2 \ge n_3 \ge\ldots > 0$ such that $\sum n_i=n$. We write $\pi=(n_1,n_2,\ldots)$, and we use the notation $|\pi|$ for the number of elements in $S_n$ which have $\pi$ as partition-type. Moreover, we let $\|\pi\|=n!/|\pi|$, so if $\sigma\in S_n$ has partition-type $\pi$, and if $Z(\sigma)$ is its centralizer, then $\|\pi\|$ is the order of $Z(\sigma)$.

To a partition we may associate the element $\l_\pi=\l^{n_1}\l^{n_2}\ldots$, as well as $\b_\pi$ and $\Psi_\pi$ defined similarly. 

Since we can express the $\b^i$'s in terms of the $\P^k$'s, it follows that the monomials $\P_\pi$, for various $\pi$, generate $\L\otimes\Q$. For dimension reasons (compare the ranks using theorem \ref{thm:lambdapoly}), we conclude that $$\L\otimes\Q=\Q[\P^1,\P^2,\ldots,\P^k,\ldots].$$

Now, consider an element $\b_H\in\L$, that is consider the permutation representation associated to the $S_n$-set $S_n/H$, for a certain subgroup $H$ of $S_n$. From what we have just said, it must be possible to write $\b_H$ as a linear combination with rational coefficients of the operations $\P_\pi$. The precise expression is (\cite{knut}, proposition on p145 and theorem on p139):
$$\b_H=\sum_\pi \frac{1}{\|\pi\|} \chi_{S_n/H}(\pi) \P_\pi.$$

Here $\chi_{S_n/H}(\pi)$ stands for the value of the character of the representation under consideration on the conjugacy class of elements which have partition-type $\pi$. 

It is possible to take this last equation, or rather the collection of equations obtained by letting $H$ vary over all subgroups of $S_n$, as a definition of the elements $\P_\pi$ for $\pi$ a partition of $S_n$. If one proceeds thus, one must work to show that each $\P^k$ is in $\L$ and not just in $\L\otimes\Q$. This is precisely what has been done to define Adams operations in a $\b$-ring.

\subsection{$\P$-operations in a $\b$-ring.}

To start with, the ring $\B$ possesses elements $\b^i$, as we have explained, and we may thus define elements $\P^k$ and even $\P_\pi\in \B$ using the same device of the log-derivative exactly as we have done in $\L$. These are not to be mistaken with the operations $\P_H$, for $H$ a subgroup of $S_n$, to be defined presently.

We shall need the following piece of notation. For $H$ a subgroup of $S_n$, we write $\|H\|$ for the order of the normalizer of $H$ in $S_n$. Moreover, $\phi_{S_n/H}$ denotes the supercharacter of the $S_n$-set $S_n/H$.

Then one has (\cite{morrisadams}):

\begin{thm} The following system of equations in $\B\otimes \Q$:
$$\b_H = \sum_K \frac{1}{\|K\|} \phi_{S_n/H}(K) \Psi_K,$$ where $K$ runs through a set of representatives for the conjugacy classes of subgroups of $S_n$, has a unique set of solutions $\Psi_K \in \B$.
\end{thm}

We call the elements $\P_K$, as well as the elements $\P_\pi$, the Adams operations in $\B$. They operate on any $\b$-ring.

\begin{rmk} Morris and Wensley in {\em loc cit} are not stating the result in this form, though it follows easily from their explicit formulae: indeed they can express $\P_K$ in terms of the $\b_H$ with the use of Moebius functions. Note also that they define elements $\psi_K$ and we have $\P_K=\| K \| \psi_K$. Let us caution the interested reader that there is a certain confusion in this article between "torsion-free ring", "ring of characteristic $0$", and "$\Q$-algebra"; moreover it is never quite clear whether the computations are done in $\B$ or $\B\otimes\Q$.
\end{rmk}

To a subgroup $H$ of $S_n$, we may associate a partition of $n$ by considering the orders of the orbits of $H$ in the set $\{1;\ldots;n\}$. We call this the {\em partition-type} of $H$ and we write $\pi_H$.

The first part of the following proposition says that $\P_\pi$ is a sort of average of the $\P_H$ for those subgroups $H$ of $S_n$ whose partition-type is $\pi$. The second part asserts that if $R$ is a $\l$-ring viewed as a $\b$-ring, then we can express all of the $\P_H$ operations on $R$ in terms of the $\P_\pi$'s.

\begin{prop}\label{prop:adams} The various Adams operations in $\B$ are related as follows.
\begin{enumerate}
\item We have the equality in $\B\otimes\Q$:
$$\P_\pi=\sum_{\pi_H=\pi} \frac{\| \pi \|}{\| H \| } \P_H$$ summing over representatives of the conjugacy classes of subgroups having $\pi$ as partition-type.

\item If $R$ is a $\l$-ring, then the operations $\P_H$ in the $\b$-ring $\b(R)$ may be described thus: When $H$ is not cyclic, $\P_H=0$. If $C$ is cyclic with partition-type $\pi$, then $\P_C=\frac{\| C \|}{\| \pi \|} \P_\pi$.

\end{enumerate}
\end{prop}

\subsection{Comments.}

The study of Adams operations on $\l$-rings has had spectacular applications, such as the solution by Adams himself of the vector field problem on spheres \cite{adams}. The Adams conjecture, now a celebrated theorem by Quillen and Sullivan and others, is stated purely in terms of operations $\P^k$ and has led to the successful computation of the image of $J$.

It can be said with confidence that the behavior of the Adams operations is the most important thing that one ought to know about a given $\l$-ring. By analogy, it is expected that the operations $\P^k$ (or $\P_\pi$) and $\P_K$ hold much information about $\b$-rings.

In the next section we shall embark on a cohomological study of the operations in the ring $\B$. The elements that we are able to deal with are the Adams operations $\P^k$, and as a result we know the behavior of $\P_\pi$ in cohomology. We hope that the reader is convinced, in particular with the help of proposition \ref{prop:adams}, that this will tell us something significant.

Before we begin, let us have a look at homotopy groups.

\subsection{The image of $J$.}

The "image of $J$" is a graded group sitting in the stable homotopy groups of spheres. Let us recall briefly the classical definition of $J$.

Pick $f \in SO_n$. Compactifying at infinity, we may think of $f$ as a pointed map $f: S^n \to S^n$. Thus if $\alpha \in \pi_k(SO_n)$, we may consider the adjoint map $S^k \wedge S^n = S^{k+n} \to S^n$ built from $f$ and $\alpha$. We obtain a homomorphism $\pi_k(SO_n)\to \pi_{n+k}(S^n)$. This passes to the limit over $n$, and we end up with:
$$J: \pi_k(SO) \to \pi_k^s(S^0).$$

The image of the homomorphism $J$ has been completely determined, though we will only state a partial result (essentially forgetting the $2$-primary part). The following theorem is due to Mitchell \cite{miche}.

\begin{thm} Let $p$ be an odd prime, and let $q$ be another prime generating the group of units in $\Z/(p^2)$. The natural map
$\bs \to BGL_\infty(\F_q)^+$ has a section at $p$; that is, the map between localised spaces $$(\bs)_{(p)} \to BGL_\infty(\F_q)^+_{(p)}$$ has a section. Moreover, on homotopy groups, this map induces an isomorphism from $(Im J)_{(p)}$ onto $\pi_*(BGL_\infty(\F_q))_{(p)}$.
\end{thm}

Thus for an odd prime, $(Im J)_{(p)}$ is zero in even degrees, and is isomorphic to $\Z/(q^i - 1)\otimes \Z_{(p)}$ in degree $2i - 1$, see example \ref{ex:adamsfinfield}. Moreover, it is a direct summand in $(\pi_*^s)_{(p)}$. Thus we obtain the following immediately (\ref{ex:adamsfinfield} and \ref{prop:ktheory}):

\begin{thm}\label{thm:imj} Let $p$ be an odd prime. Let $\P^k:\bs \to \bs$ be the map representing the Adams operation on $\tilde{\pi}^0(-)$, and let 
$$\P^k_* : (Im J)_{(p)} \to (Im J)_{(p)}$$ be the map induced on the image of $J$ using the projection $(\pi_*^s)_{(p)} \to (Im J)_{(p)}$. Then $\P^k_*$ is multiplication by $k^i$ in degree $2i -1$ (and zero otherwise).
\end{thm}

\section{Effect in cohomology}

The operations $\b_H$ on $\tilde{\pi}^0(-)$ are represented by maps
$$\b_H: \bs\to\bs$$ In this section we intend to give some information on the induced maps $\b_H^*$ on $H^*(\bs,k)$ where $k=\F_2$ (it is expected that the cases $k=\Z$ or $k=\F_p$ will give a similar answer); more precisely we shall consider the Adams operations.

We need to motivate the form that the result will take, and for this we recall what is known in the case of the operations $\l^n$ in $\L$ which act on $\tilde{K}^0(-)$ and are represented by maps $\l^n:BU\to BU$.

\subsection{The case of $K$-theory.}

Let $c\in H^*(BU,\Z)$. We may think of $c$ as a characteristic class for representations of finite groups: given a finite group $G$ and a reresentation $\rho: G \to U(n)$, compose $B\rho$ with the natural map $BU(n)\to BU$ to obtain $f: BG\to BU$. Then $f^*(c)\in H^*(BG,\Z)$ is denoted by $c(\rho)$. In fact, it is well-known that a class such as $c$ is completely determined by the values $c(\rho)$ for various $G$ and $\rho$.

We may extend $c$ to a characteristic class for virtual representations, as follows. The cohomology of $BU$ is known:
$$H^*(BU,\Z)=\Z[c_1,c_2,\ldots,c_n,\ldots].$$ The classes $c_i$ are known as Chern classes. Consider then the sum $$c_t=1 + c_1t + c_2t^2 + \ldots \in 1 + tH^*(BU,\Z)[[t]].$$ More precisely, if we give $t$ the weight $-2$, then $c_t$ belongs to the subset of homogeneous elements of weight $0$ in  $1 + tH^*(BU,\Z)[[t]]$, which is a multiplicative group that we shall denote by $H^\times(BU,\Z)$. We can see $c_t$ again as a characteristic class for representations of $G$ with values in $H^\times(BG,\Z)$. We call $c_t$ an {\em exponential} characteristic class since it enjoys the following properties:
\begin{itemize}
\item $c_t(\textrm{trivial representation})=1$,
\item $c_t(E\oplus F)=c_t(E)c_t(F)$.
\end{itemize}

From this, it follows easily that $c_t$ extends to a homomorphism of groups:
$$c_t: R(G) \to H^\times (BG,\Z).$$ Thus, any $c\in H^*(BU,\Z)$, being a polynomial in the Chern classes, extends to $R(G)$ when viewed as a characteristic class.

Of particular interest to us are the elements $\chi^i$ defined by
$$\chi^i=P(c_1,c_2,\ldots,c_i)$$ where $P$ is that polynomial such that
$$\P^i=P(\l^1,\l^2,\ldots,\l^i).$$ Equivalently, the $\chi^i$'s are the coefficients of the power series $-t\frac{c_{-t}'}{c_{-t}}.$ It follows exactly as in the case of the $\P$-operations that the rational cohomology of $BU$ is a polynomial ring on the classes $\chi^i$. Moreover, these may be put together into a map:
$$\chi: R(G) \to H^*(BG,\Z)$$ by taking $\chi(\rho)=\sum_i \chi^i(\rho)$. This map $\chi$ is a homomorphism of groups: this comes from the fact that $c_t$ is exponential, and is proved exactly as the additivity of the Adams operations.

Now, since $H^*(BU,\Z)$ is torsion-free, the effect of a map $f:BU\to BU$ on cohomology may be read from the values $f^*(\chi^i)$. For the same reason, since any element of $\L$ is a polynomial in the Adams operations with rational coefficients, it is sufficient to know the effect of the maps $(\P^k)^*$ in order to determine completely the behavior of $\alpha^*$ for any $\alpha\in\L$.

So in the end all we need to know is the value of $(\P^k)^*(\chi^i)$, and from the point of view of characteristic classes, what we require is a formula for $\chi^i(\P^k(\rho))$ in terms of the classes $\chi^j(\rho)$. This result is due to Adams. Let us consider $H^*(BG,\Z)$ as a $\P$-ring as in example \ref{ex:psigraded}. Then for any representation $\rho$:
$$\chi^i(\P^k(\rho))=\P^k(\chi^i(\rho)).$$ In other words one has
$$(\P^k)^*(\chi^i)=\P^k(\chi^i) \qquad (=k^{2i}\chi^i)$$ in the integral cohomology of $BU$.

\begin{rmk} In the literature, one often considers the classes $$Ch^i=\frac{1}{i!}\chi^i\in H^*(BU,\Q)$$ and the corresponding map $Ch$ which is called the {\em Chern character}. This has the advantage of being a multiplicative map, and moreover when $X$ is compact and connected $Ch$ yields an isomorphism $K^0(X)\otimes\Q \to H^*(X,\Q)$. While these properties are very useful in other contexts, we prefer to consider the classes without denominators which live in the integral cohomology. Indeed, in the sequel we shall be dealing with $\bs$ whose rational cohomology is zero.
\end{rmk}

\subsection{Effect of the $\b$-operations in cohomology.}

Let $c\in H^*(\bs,k)$. We may think of $c$ as a characteristic class for set-theoretic representations of finite groups: given a finite group $G$ and a reresentation $\rho: G \to S_n$, compose $B\rho$ with the natural map $BS_n\to \bs$ to obtain $f: BG\to \bs$. Then $f^*(c)\in H^*(BG,k)$ is denoted by $c(\rho)$. In fact, it is well-known that a class such as $c$ is completely determined by the values $c(\rho)$ for various $G$ and $\rho$ (we shall give more details on this in the course of the proof of theorem \ref{thm:psicoho} below).

Suppose given a sequence $c_1, c_2, \ldots, c_n, \ldots$ of elements of $H^*(\bs,k)$, with $|c_i|=2i$ (or $|c_i|=i$ when $k=\F_2$). Consider then the sum $$c_t=1 + c_1t + c_2t^2 + \ldots \in 1 + tH^*(\bs,k)[[t]].$$ More precisely, if we give $t$ the weight $-2$ (or $-1$ when $k=\F_2$), then $c_t$ belongs to the subset of homogeneous elements of weight $0$ in  $1 + tH^*(\bs,k)[[t]]$, which is a multiplicative group that we shall denote by $H^\times(\bs,k)$. We can see $c_t$ again as a characteristic class for representations of $G$ with values in $H^\times(BG,k)$. We call $c_t$ an {\em exponential} characteristic class when it enjoys the following properties:
\begin{itemize}
\item $c_t(\textrm{trivial representation})=1$,
\item $c_t(E\coprod F)=c_t(E)c_t(F)$.
\end{itemize}

(See \ref{ex:segalstret} below for examples). It follows easily that an exponential class $c_t$ extends to a homomorphism of groups:
$$c_t: A(G) \to H^\times (BG,k).$$ 

We shall be particularly interested in the elements $\chi_c^i$ given by the coefficients of the power series$-t\frac{c_{-t}'}{c_{-t}}.$ These may be put together into a map:
$$\chi_c: A(G) \to H^*(BG,\Z)$$ by taking $\chi_c(\rho)=\sum_i \chi_c^i(\rho)$. This map $\chi_c$ is a homomorphism of groups which we call the {\em character of $c_t$}.

\begin{thm}\label{thm:psicoho}
Let $k=\F_2$, let $c_t\in H^\times(\bs,\F_2)$ be an exponential characteristic class, and let $\chi_c$ denote the character of $c_t$. Then $$(\P^k)^*(\chi_c^i)=\P^k(\chi_c^i)=k^i\chi_c^i.$$ In other words, for any finite group $G$ and any $\rho\in A(G)$, we have
$$\chi_c^i(\P^k(\rho))=\P^k(\chi_c^i(\rho)).$$ 
\end{thm}

This result will certainly have to be refined in the cases $k=\F_p$ and $k=\Z$, though we hope that the answer will be similar in spirit. Of course when we work in mod $2$ cohomology, what we need to prove reduces to:
\begin{itemize}
\item $\chi_c^i(\P^k(\rho))=0$ when $k$ is even,
\item $\chi_c^i(\P^k(\rho))=\chi_c^i(\rho)$ when $k$ is odd. 
\end{itemize}

\begin{proof} Let us consider a particular case first. Let $A_n=(\Z/2)^n$ be elementary abelian, so that the order of $A_n$ is $2^n$. For any finite group $G$ and any $x\in A(G)$, we have 
$\P^k(x)=\P^d(x)$ where $d=gcd(k,|G|)$ (\cite{gay}, corollary 4.3). So if $k$ is odd and $\rho\in A(A_n)$, we have $\P^k(\rho)=\P^1(\rho)=\rho$ and we have indeed $\chi_c^i(\P^k(\rho))=\chi_c^i(\rho)$.

Now let $k$ be even, and $d=gcd(k,|A_n|)=2^r$ with $r\ge 1$. By the same result as above we have $\P^k(\rho)=\P^{2^r}(\rho)$. However, if $P$ denotes that polynomial such that $\P^i=P(\l^1,\ldots,\l^i)$, then $P(\epsilon_1,\ldots,\epsilon_i)=t_1^i + t_2^i + \ldots + t_i^i$ when $\epsilon_j$ is the $j$-th symmetric function in the variables $t_1,\ldots,t_i$. In particular, $\P^{2^r}=(\l^1)^{2^r} mod ~2$ since $t_1^{2^r} + t_2^{2^r} + \ldots + t_{2^r}^{2^r}=(t_1 + \ldots t_{2^r})^{2^r} mod ~2$. So $\P^k(\rho)=\P^{2^r}(\rho)=\rho^{2^r} mod~ 2$. Since $\chi_c^i$ is additive and takes values in mod $2$ cohomology, it is enough to show $\chi_c^i(\rho^{2^r})=0$. By additivity again (of $\chi_c^i$ and $\P^k$), it is enough to consider the case $\rho=A_n/U$ for a subgroup $U$ of $A_n$. We assume that the index of $U$ in $A_n$ is $>1$ since any character is $0$ on trivial representations anyway.

Consider then the $U$-set $X=(A_n/U)_{|U}^{2^r - 1}$. Since $A_n$ is abelian, $X$ is a trivial $U$-set, so it splits as $2^s$ copies of the one-point trivial $U$-set. If we consider the {\em transfer} $X^{U\to A_n}$, then we see that it splits into $2^s$ copies of the same representation (namely $A_n/U$). Therefore $\chi_c^i(X^{U\to A_n})=0$.

However, $X^{U\to A_n}$ is none other than $(A_n/U)^{2^r}$. To see this, use the equality $$(M_{|U} \times N)^{U\to A_n}=M \times N^{U\to A_n}.$$ The latter may be proved directly, or see \cite{sieben}, 2.19.2. Using this first with $M=(A_n/U)^{2^r - 2}$ and $N=(A_n/U)_{|U}$ and then with $M=A_n/U$ and $N=U/U$, we do obtain the identification of $(A_n/U)^{2^r}$ with $X^{U\to A_n}$. Therefore $\chi_c^i((A_n/U)^{2^r})=0$, and the theorem holds when $G=A_n$.

To deal with the general case, let us recall a classical result. The homology $H_*(\bs,\F_2)=H_*(BS_\infty,\F_2)$ can be turned into a ring by using the natural maps
$$S_n\times S_m \to S_{n+m}.$$ Moreover, it is known that $H_*(BS_\infty,\F_2)$ is generated as a ring by the elements $\iota^n_*(x)$ where $x\in H_*(BA_n,\F_2)$ and $\iota^n : A_n \to S_{2^n}$ is the regular representation of $A_n$ (we simplify the notations and write $\iota^n_*$ instead of $B\iota^n_*$). Consider then $$\rho=\iota^n\times \iota^m: A_{n+m}=A_n\times A_n \to S_{2^n}\times S_{2^m} \to S_{2^n+2^m} \to S_\infty.$$ The product $\iota^n_*(x)\cdot \iota^m_*(y)$ in $H_*(BS_\infty,\F_2)$ is by definition $\rho_*(x\otimes y)$.

Using this repeatedly, we see that $H_*(BS_\infty,\F_2)$ is generated additively by elements of the form $\rho_*(x)$ where $x\in H_*(BA_N,\F_2)$ and $\rho: BA_N \to BS_\infty$ comes from a representation of $A_N$. Dually, a cohomology class $c\in H^*(\bs,\F_2)$ is completely determined by the values
$$\<\rho_*(x),c\>=\<x,\rho^*(c)\>=\<x,c(\rho)\>,$$ that is, $c$ is determined by the cohomology classes $c(\rho)$. From the example we have just considered, it follows that $(\P^k)^*(\chi_c^i)$ and $\P^k(\chi_c^i)=k^i\chi_c^i$ must be equal.
\end{proof}

\begin{ex}\label{ex:segalstret} Segal and Stretch \cite{segalstret} have constructed an infinite family $c_t^{(k)}$ of exponential characteristic classes. These are finer than Stiefel-Whitney or Chern classes (which are also exponential), in that they can distinguish between representations in $A(G)$ which have the same image in $R(G)$.
\end{ex}

\bibliography{myrefs}
\bibliographystyle{siam}
\end{document}